%% file: ad_lw__arXiv.tex
\documentclass[]{scrartcl}

\usepackage[utf8]{inputenc}  
\usepackage[USenglish]{babel}
\usepackage{csquotes}

\usepackage[svgnames]{xcolor}

\usepackage[a4paper,top=22mm,bottom=20mm,inner=25mm,outer=20mm]{geometry}

\usepackage[%
  backend=bibtex,bibencoding=ascii,
  style=numeric-comp,
  giveninits=true, uniquename=init, 
  natbib=true,
  url=true,
  doi=true,
  isbn=false,
  backref=false,
  maxnames=99,
  ]{biblatex}
\usepackage{amsmath}
\allowdisplaybreaks
\numberwithin{equation}{section}
\usepackage{amssymb}
\usepackage{commath}
\usepackage{mathtools}
\usepackage{bbm}
\usepackage{nicefrac}
\usepackage{subdepth}
\usepackage[frozencache,cachedir=minted-cache]{minted}
\usepackage{minted}

\usepackage{algorithm}
\usepackage{algpseudocode}

\usepackage{amsthm}
\usepackage{thmtools}
\usepackage{etoolbox}
\makeatletter
\patchcmd{\thmt@setheadstyle}
 {\bgroup\thmt@space}
 {\thmt@space}
 {}{}
\patchcmd{\thmt@setheadstyle}
 {\egroup\fi}
 {\fi}
 {}{}
\makeatother
\declaretheoremstyle[
  bodyfont=\normalfont\itshape,
  headformat=\NAME\ \NUMBER\NOTE,
]{myplain}
\declaretheoremstyle[
  headformat=\NAME\ \NUMBER\NOTE,
]{mydefinition}
\newcommand{\envqed}{{\lower-0.3ex\hbox{$\triangleleft$}}}
\declaretheorem[style=myplain,numberwithin=section]{theorem}

\declaretheorem[style=mydefinition,numberlike=theorem,qed=\envqed]{definition}

\usepackage[plainpages=false,pdfpagelabels,hidelinks,unicode]{hyperref}

\usepackage{color}
\usepackage{graphicx}
\usepackage[small]{caption}
\usepackage{subcaption}

\begingroup\expandafter\expandafter\expandafter\endgroup
\expandafter\ifx\csname pdfsuppresswarningpagegroup\endcsname\relax
\else
\fi

\usepackage{booktabs}
\usepackage{rotating}
\usepackage{multirow}

\usepackage{enumitem}

\usepackage[T1]{fontenc}

\usepackage{newpxtext,newpxmath}

\usepackage{lineno}

\input{macros}

\newcommand{\orcid}[1]{ORCID:~\href{https://orcid.org/#1}{#1}}
\usepackage{authblk}

\newenvironment{keywords}{\par\textbf{Key words.}}{\par}

\title{Automatic differentiation for \mrev{performing the Cauchy-Kovalevskaya procedure in} Lax-Wendroff type discretizations}

\author[1]{Arpit~Babbar\thanks{\orcid{0000-0002-9453-370X}}}
\affil[1]{Institute of Mathematics, Johannes Gutenberg University Mainz, Staudingerweg 9, 55128 Mainz, Germany}

\author[1,2]{Valentin~Churavy\thanks{\orcid{0000-0002-9033-165X}}}

\author[2]{Michael~Schlottke-Lakemper\thanks{\orcid{0000-0002-3195-2536}}}
\affil[2]{High-Performance Scientific Computing, Centre for Advanced Analytics and Predictive Sciences, University of Augsburg, Germany}

\author[1]{Hendrik~Ranocha\thanks{\orcid{0000-0002-3456-2277}}}

\makeatletter
\hypersetup{pdfauthor={Arpit~Babbar, Valentin~Churavy, Michael~Schlottke-Lakemper, Hendrik~Ranocha}} 
\hypersetup{pdftitle={Automatic differentiation for
\mrev{performing the Cauchy-Kovalevskaya procedure in} Lax-Wendroff type
discretizations}} 
\makeatother

\begin{document}

\maketitle

\begin{abstract}
\noindent
  \input{abstract.tex}
\end{abstract}



\section{Introduction}

\mrev{This work is a contribution to }Lax-Wendroff (LW) methods \mrev{which are based on the Cauchy-Kovalevskaya procedure to} obtain order $N+1$ accuracy for solving the conservation law $$\uu_t + \pf (\uu)_{x} = \bzero$$ by approximating the truncated time averaged flux
\mrev{\begin{equation}
\F = \sum_{m \subindex 0}^N \frac{\mrev{\Delta t^m}}{(m + 1) !} \partial_t^m  \pf.
\label{eq:tavgflux.intro}
\end{equation}}
\mrev{The first Lax-Wendroff method was the second order finite difference method introduced in~\cite{Lax1960} for solving hyperbolic conservation laws.
The method was extended to higher order accuracy by combining the Taylor expansion in time with the Cauchy-Kovalevskaya procedure that replaces temporal derivatives in the time averaged flux~\eqref{eq:tavgflux.intro} by spatial derivatives using the conservation law.
This procedure was first introduced in~\cite{Qiu2003,Qiu2005b} where the obtained spatial derivatives were discretized using WENO and discontinuous Galerkin (DG) methods respectively.
This is also considered to be a particular case of ADER (Arbitrary high order schemes
using DERivatives) methods (see~\cite{Titarev2002}), although the latter term is often also used to refer to methods that use a local space-time predictor based on a weak formulation~\cite{dumbser2008,dumbser2008a}.
Both ADER and LW methods have been studied extensively in the literature as high-order one-step time discretizations for hyperbolic conservation laws.
Such methods are an important alternative to the widely used Runge-Kutta time discretizations which require multiple stages to achieve high-order accuracy in time, requiring inter-element communication at each stage.
}
High-order LW methods were introduced in~\cite{Qiu2003,Qiu2005b} based on the chain rule \mrev{with analytically computed derivatives}.
For example, the first and second derivatives are successively obtained by
$$\pf_t = \pf'\uu_t, \quad \uu_t = - \pf_x$$
and
$$\pf_{t  t} = \left( \pf'  \uu_t \right)_t = \pf''  \uu_t \uu_t + \pf'  \uu_{t  t}, \quad \uu_{t t} = - \left( \pf_t \right)_x.$$
The \textit{approximate LW procedure} was proposed in~\cite{Zorio2017} as an alternative using finite difference (FD) approximations for computing $\partial_t^m \pf$.
The approximate LW procedure was extended to discontinuous Galerkin (DG) methods in~\cite{Burger2017} and was demonstrated to be faster than the LW DG methods of~\cite{Qiu2005b}.
The reason for the observed speed-up was given to be the computation of the analytic flux Jacobians in the application of chain rule, which are not needed in an FD approximation.
\mrev{In this work, we use the Lax-Wendroff Flux Reconstruction (LWFR) method of~\cite{babbar2022lax} which uses the flux reconstruction (FR) spatial discretization~\cite{Huynh2007}.
However, the work presented here is also directly applicable to DG methods due to the equivalence of FR and DG methods shown in~\cite{Huynh2007}.}
We propose using automatic differentiation (AD) for computing the temporal derivatives $\partial_t^m  \pf$~\eqref{eq:tavgflux.intro}.
This is equivalent to the application of the chain rule used in~\cite{Qiu2003,Qiu2005b} but does not require analytic Jacobians as it computes directional derivatives and applies to any choice of the flux $\pf$.
Thus, as observed in our numerical experiments, the performance of AD is always \mrev{on par or slightly faster than} that of the approximate LW procedure.
The implementation of the AD approach is identical for any desired order of accuracy, unlike the approximate LW procedure, where different FD methods are required for higher order of accuracy.
Another benefit of AD over the approximate LW procedure is that, for the latter, the admissibility preserving procedure of~\cite{babbar2024admissibility,babbar2025generalized} is insufficient for conservation laws whose fluxes cannot be computed on physically inadmissible states.
This was first observed in~\cite{basak2025} for relativistic hydrodynamic equations, where additional positivity correction was performed for the local time averaged flux computation (predictor step) to be able to use the approximate LW procedure.
Since our proposed AD approach only uses flux evaluations of the approximate solution which is physically admissible following the procedure of~\cite{babbar2024admissibility}, it does not need positivity correction in the predictor step.

The rest of the paper is organized as follows.
\mrev{Section~\ref{sec:lwfr} reviews the LW method with flux reconstruction spatial discretization and the admissibility preservation procedure of~\cite{babbar2024admissibility}.}
In Section~\ref{sec:lw}, we review the approximate LW procedure and introduce the AD approach.
In Section~\ref{sec:results}, we present numerical results to demonstrate the accuracy and admissibility preservation of the AD approach and compare its performance with the approximate LW procedure.
\mrev{Conclusions are drawn from the paper in Section~\ref{sec:conclusions}.
Some background on automatic differentiation along with a basic implementation using dual numbers is given in Appendix~\ref{sec:ad.appendix}.}

\section{\mrev{Lax-Wendroff flux reconstruction with subcell based admissibility preservation}} \label{sec:lwfr}

The Lax-Wendroff flux reconstruction (LWFR) method of~\cite{babbar2022lax} is reviewed in this section along with the admissibility preservation procedure of~\cite{babbar2024admissibility,babbar2025generalized}.
We consider the one-dimensional hyperbolic conservation law
\begin{equation}
\uu_t + \pf (\uu)_{x} = \bzero. \label{eq:con.law}
\end{equation}
A solution $\uu$ is said to be physically admissible if it belongs to the set $\Uad$ defined by
\begin{equation}
	\label{eq:uad.form} \Uad = \{ \uu \in \re^p : \ad_k (\uu) > 0, 1 \le k \le	K\},
\end{equation}
where $\{\ad_k\}$ depend on the equations. For Euler's equations, $K = 2$ and $\ad_1, \ad_2$ are density, pressure
functions respectively.
We will divide the computational domain $\Omega$ into disjoint elements $\{\Omega_e\}$, with $\Omega_e = [x_{\emh}, x_{\eph}]$ so that $\Delta x_e = x_{\eph} - x_{\emh}$.
Then, we map each element $\Omega_e$ to the reference element $[0, 1]$ by $x \mapsto \frac{x - x_{e-1/2}}{\Delta x_e} =: \xi$.
Inside each element, we approximate the solution by degree $N \ge 0$ polynomials belonging to the set $\poly_N$.
For this, we choose $N + 1$ distinct nodes $\{ \xi_i \}_{i \subindex 0}^N$ in $[0,1]$ which will be taken to be Gauss-Legendre (GL) points in this work, and will also be referred to as \textit{solution points}.
The solution inside an element $\Omega_e$ is given by $\uu_h (\xi, t) = \sum_{p = 0}^N \uep (t) \ell_p(\xi)$ where $\{\ell_p\}$ are Lagrange polynomials of degree $N$ defined to satisfy $\ell_p (\xi_q) = \delta_{pq}$.
Note that the coefficients $\uep$ which are the basic unknowns or \textit{degrees of freedom} (dof), are the solution values at the solution points.
Following equation~\eqref{eq:lwtay}, we need to specify the construction of the time averaged flux~\eqref{eq:tavgflux}.
The first step of approximating~\eqref{eq:lwtay} is the predictor step where a locally degree $N$ approximation $\F^\delta$ of the time averaged flux is computed by the Cauchy-Kovalevskaya procedure which we now review.
\mrev{Using a Taylor expansion in time around $t = t_n$, we can write the solution at the next time level as $$\uu^{n + 1} = \uu^n + \sum_{m = 1}^{N + 1} \frac{\Delta t^m}{m!} \partial_t^m  \uu^n + O (\Delta t^{N + 2}),$$ where $N$ is the polynomial degree of the DG or flux reconstruction (FR) method.
Since the spatial error is expected to be $O (\Delta x^{N + 1})$, we retain terms up to $O (\Delta t^{N + 1})$ in the Taylor expansion, so that the overall formal accuracy is of order $N + 1$ in both space and time.
Using the conservation law, we re-write temporal derivatives of
the solution in terms of spatial derivatives of the flux as $\partial_t^m  \uu = - (\partial_t^{m - 1}  \pf)_x$, so that
\begin{align}
\uu^{n + 1} & = \uu^n - \sum_{m = 1}^{N + 1} \frac{\Delta t^m}{m!}
(\partial_t^{m - 1}  \pf)_x + O (\Delta t^{N + 2}) =  \uu^n - \Delta t \pd{\F}{x} (\uu^n) + O (\Delta t^{N + 2}), \label{eq:lwtay} \\
\F & = \sum_{m = 0}^N \frac{\Delta t^m}{(m + 1) !}
\partial_t^m  \pf = \pf + \frac{\Delta t}{2} \partial_t  \pf
+ \ldots + \frac{\Delta t^N}{(N + 1) !} \partial_t^N  \pf.
\label{eq:tavgflux}
\end{align}
Note that $\F = \F({\uu^n})$~\eqref{eq:tavgflux} is an approximation to the time average flux in the interval $[t_n, t_{n + 1}]$~\cite{babbar2022lax}.}
In the original LWFR method of~\cite{babbar2022lax}, the approximate LW procedure of~\cite{Zorio2017,Burger2017} was used for approximating the temporal derivatives of the flux in the Cauchy-Kovalevskaya procedure~\eqref{eq:tavgflux}.
The review of the approximate LW procedure and the proposed AD approach is given in Section~\ref{sec:lw}.
Then, as in the standard RKFR scheme, we perform the flux reconstruction procedure on $\F^\delta$ to construct a locally degree $N+1$ and globally continuous flux approximation $\F_h (\xi)$.
After computing $\F_h$, truncating equation~\eqref{eq:lwtay}, the solution at the nodes is updated by a collocation scheme as follows
\begin{equation}
\uep^{n + 1} = \uep^n - \frac{\Delta t}{\Delta x_e}  \od{\F_h}{\xi} (\xi_p), \qquad 0 \le p \le N. \label{eq:uplwfr}
\end{equation}
This is the single stage Lax-Wendroff update scheme for any order of accuracy.
In order to preserve admissibility of the solution, we use the subcell based admissibility preservation procedure of~\cite{babbar2024admissibility,babbar2025generalized}.
In~\cite{babbar2024admissibility}, it was used in combination with a subcell based blending limiter inspired from~\cite{hennemann2021} that is used for controlling spurious oscillations.
Here we describe the generalized admissibility preservation procedure introduced in~\cite{babbar2025generalized} which can be done independently of the choice of limiter.

The admissibility preserving property of the conservation law, also known as convex set preservation property since $\Uad$ is convex, can be written as
\begin{equation}
	\label{eq:conv.pres.con.law} \uu (\cdummy, t_0) \in \Uad \qquad
	\Longrightarrow \qquad \uu (\cdummy, t) \in \Uad, \qquad t > t_0.
\end{equation}
Thus, an admissibility preserving FR scheme is one which preserves admissibility at all solution points. It is defined precisely as follows.
\begin{definition}\label{defn:admissibility.preserving}
A flux reconstruction scheme is said to be admissibility-preserving if
\[
\uu_{e,p}^n \in \Uad \quad \forall e,p \qquad \implies \qquad \uu_{e,p}^{n+1} \in \Uad \quad \forall e,p,
\]
where $\Uad$ is the admissible set~\eqref{eq:uad.form} of the physical equation~(\ref{eq:uad.form}).
\end{definition}
In this work, as in~\cite{babbar2024admissibility}, we study the admissibility preservation in means property of the LWFR scheme defined as
\begin{definition}\label{defn:admissibility.preserving.means}
A flux reconstruction scheme is said to be admissibility-preserving in means if
\begin{equation}
\uu_{e,p}^n \in \Uad \quad \forall e,p \qquad \implies \qquad \au_{e}^{n+1} \in \Uad \quad \forall e,
\end{equation}
where $\au_{e}$ is the element mean of $\uu_h$ in element $e$ computed as
\[
\au_e = \sum_{p=0}^N w_p \uu_{e,p},
\]
where $w_p$ are the quadrature weights associated with the solution points, and $\Uad$ is the admissible set~\eqref{eq:uad.form} of the equations of interest~(\ref{eq:uad.form}).
\end{definition}
Once the scheme is admissibility preserving in means, the scaling limiter of~{\cite{zhang2010c}} can be used to obtain an admissibility preserving scheme.
The following property of the LWFR scheme will be crucial in obtaining the admissibility preservation in means
\begin{equation}
\label{eq:fravgup} \au_e^{n + 1} = \au_e^n - \frac{\Delta t}{\Delta x_e}  (\F_{\eph} - \F_{\emh}).
\end{equation}
\subsection{\mrev{Limiting time average flux}} \label{sec:flux.correction}
As in~\cite{babbar2025generalized}, we define \textit{fictitious finite volume updates}
\begin{equation}
	\begin{split}
		\atu_{e, 0}^{n + 1} & = \uez^n
    - \frac{\Delta t}{w_0 \Delta x_e}
		[\pf_{\half}^e - \F_{\emh}^{\text{LW}}],\\
		\atu_{e, p}^{n + 1} & = \uep^n - \frac{\Delta t}{w_p \Delta x_e}
		[\pf_{\pph}^e - \pf_{\pmh}^e], \qquad 1 \le p \le N - 1,\\
		\atu_{e, N}^{n + 1} & = \ueN^n - \frac{\Delta t}{w_N \Delta x_e}
		[\F_{\eph}^{\text{LW}} - \pf_{\Nmh}^e],
	\end{split} \label{eq:low.order.update}
\end{equation}
where $\pf_{\pph}^e = \pf ( \uep^n, \ueppone^n )$ is an admissibility preserving finite volume numerical flux.
Then, note that
\begin{equation}\label{eq:cell.avg.decomp}
	\avg{\uu}_e^{n + 1} = \sum_{p = 0}^N w_p  \atu_{e, p}^{n + 1}.
\end{equation}
Thus, if we can ensure that $\atu_{e, p}^{n + 1} \in \Uad$ for all
$p$, the scheme will be admissibility preserving in
means (Definition~\ref{defn:admissibility.preserving.means}). We do have $\atu_{e, p}^{n + 1} \in \Uad$ for $1 \le p \le N-1$ under appropriate CFL conditions because the finite volume fluxes are admissibility preserving.
In order to ensure that the updates $\atu_{e, 0}^{n + 1}, \atu_{e, N}^{n + 1}$ are also admissible, the flux limiting procedure of~{\cite{babbar2024admissibility}} is followed so that the high order numerical fluxes $\F_{e\pm\half}^{\text{LW}}$ are replaced by the \textit{blended numerical fluxes} $\F_{e \pm \half}$.
The procedure is explained here for completeness. We define an admissibility preserving lower order flux at the interface $\eph$ as
\[
\pf_{\eph} = \pf ( \uepoz^n, \ueN^n ).
\]
Note that, for an RKFR scheme using Gauss-Legendre-Lobatto (GLL) solution
points, the definition of $\atu_{e, N}^{n + 1}$ will use $\pf_{\eph}$ in place of $\F_{\eph}^{\text{LW}}$ and thus the admissibility preserving in means property (Definition~\ref{defn:admissibility.preserving.means}) will always be present.
That is the argument of~{\cite{zhang2010c}} and here we review how~\cite{babbar2025generalized} extends it to LWFR schemes by limiting $\F_{\eph}^\text{LW}$.
We will explain the procedure for limiting $\F_{\eph}^{\text{LW}}$ to obtain $\F_{\eph}$; it will be similar in the case of $\F_{\emh}$.
Note that we want $\F_{\eph}$ to be such that the following are admissible
\begin{equation}
	\begin{split}
		\atu_0^{n + 1} & = \uepoz^n - \frac{\Delta t}{w_0 \Delta x_{e + 1}}
		(\pf^{e + 1}_{\half} - \F_{\eph}),\\
		\atu_N^{n + 1} & = \ueN^n - \frac{\Delta t}{w_N \Delta x_e}  (\F_{\eph}
		- \pf^e_{\Nmh}).
	\end{split} \label{eq:low.order.tilde.update}
\end{equation}
We will exploit the admissibility preserving property of the finite volume scheme to get
\begin{align*}
\utilow_0 & = \uepoz^n - \frac{\Delta t}{w_0 \Delta x_{e + 1}}
(\pf^{e+1}_{\half} - \pf_{\eph}) \in \Uad,\\
\utilow_N & = \ueN^n - \frac{\Delta t}{w_N \Delta x_e}  (\pf_{\eph} - \pf^e_{\Nmh}) \in \Uad.
\end{align*}
Let $\{ \ad_k, 1 \le k \le K\}$ be the admissibility constraints~\eqref{eq:uad.form} of~\eqref{eq:con.law}.
The time averaged flux is limited by iterating over these constraints.
For the $k^\text{th}$ constraint, we can solve an optimization problem to find the largest $\theta \in [0,1]$ satisfying
\begin{equation}
\ad_k ( \theta \atu_p^{n + 1} + (1 - \theta)  \utilow_p ) > \epsilon_p, \qquad p = 0, N, \label{eq:optimization.problem}
\end{equation}
where $\epsilon_p$ is a tolerance, taken to be $\frac{1}{10}  \ad_k (\utilow_p)$~\cite{ramirez2021}.
The optimization problem is usually a polynomial equation in $\theta$.
If $\ad_k$ is a concave function of the conserved variables, we can follow~{\cite{babbar2024admissibility}} and use the simpler but possibly sub-optimal approach of defining
\begin{equation}
\theta = \min \left( \min_{p = 0, N} \frac{|\epsilon_p - \ad_k (\atu_p^{\text{low}, n + 1})|}{|\ad_k (\atu_p^{n + 1})  - \ad_k (\atu_p^{\text{low}, n + 1})| + \texttt{eps}} , 1 \right), \label{eq:concave.theta}
\end{equation}
where $\texttt{eps} = 10^{-13}$ is used to avoid a division by zero.
In either case, by iterating over the admissibility constraints $\{ \ad_k \}$ of the conservation law, the flux $\F_{\eph}^{\text{LW}}$ can be corrected as in Algorithm~\ref{alg:flux.correction}.
\begin{algorithm}
\caption{Flux limiting} \label{alg:flux.correction}
\begin{algorithmic}
\State $\F_{\eph} \leftarrow \F_{\eph}^{\text{LW}}$
\For{$k$=1:$K$}
\State $\epsilon_0, \epsilon_N \leftarrow \frac{1}{10}  \ad_k (\utilow_0), \frac{1}{10}  \ad_k (\utilow_N)$
\State Find $\theta$ by solving the optimization problem~\eqref{eq:optimization.problem} or by using~\eqref{eq:concave.theta} if $P_k$ is concave
\State $\F_{\eph} \leftarrow \theta \F_{\eph} + (1 - \theta) \pf_{\eph}$
\State $\atu_0^{n + 1} \leftarrow \uepoz^n - \frac{\Delta t}{w_0 \Delta x_{e + 1}}  (\pf^{e + 1}_{\half} - \F_{\eph})$
\State $\atu_N^{n + 1} \leftarrow \ueN^n - \frac{\Delta t}{w_N \Delta x_e}  (\F_{\eph} - \pf^e_{\Nmh})$
\EndFor
\end{algorithmic}
\end{algorithm}
After $K$ iterations, we will have $\ad_k (\atu_p^{n + 1}) \geq \epsilon_p$ for $p = 0, N$ and $1 \leq k \leq K$.
Once Algorithm~\ref{alg:flux.correction} is performed at all interfaces, using the numerical flux $\F_{\eph}$ in the LWFR scheme will ensure that $\atu_{e, p}^{n +1}$~\eqref{eq:low.order.update} belongs to $\Uad$ for all $p$, implying $\au_e^{n + 1} \in \Uad$ by~\eqref{eq:cell.avg.decomp}.
We can then use the scaling limiter of~\cite{zhang2010c} at all solution points and obtain an admissibility preserving LWFR scheme for conservation laws.
\section{Approximate Lax-Wendroff procedure and automatic differentiation}\label{sec:lw}

The spatial derivatives $(\partial_t^m  \pf)_x$ \mrev{in~\eqref{eq:tavgflux}} are computed by differentiating their element local polynomial approximations in the DG/FR basis~\cite{babbar2022lax}.
These are then used to obtain a local approximation of the time averaged flux $\F$\mrev{~\eqref{eq:tavgflux}} which is made global by the DG/FR procedure~\cite{babbar2022lax} and used to update the solution as in~\eqref{eq:lwtay}.
We now discuss the existing options for obtaining $\partial_t^m  \pf$ using $\{ \partial_t^k \uu\}_{k:0}^{m-1}$, and then how AD can be used as an alternative.
\mrev{The relevant background for automatic differentiation is given in Appendix~\ref{sec:ad.appendix}.}
Following~\cite{babbar2022lax}, we write
\[
\vu^{(m)} = \Delta t^m \partial_t^m  \uu, \qquad \vf^{(m)} = \Delta t^m \partial_t^m  \pf.
\]
The temporal derivatives of the flux are given by Faà di Bruno's formula~\cite{faadibruno1857}
\begin{equation} \label{eq:faa}
\vf^{(m)} = \sum_{\pi \in \Pi} \partial_{\mrev{\uu}}^{|\pi|}\pf (\uu) \prod_{B \in \pi} \vu^{(|B|)},
\end{equation}
where $\Pi$ is the set of all partitions of the set $\{1, \ldots, m\}$, $|\pi|$ is the number of blocks in the partition $\pi$, and $|B|$ is the number of elements in the block $B$. The first works using LW methods~\cite{Qiu2003,Qiu2005b} employ~\eqref{eq:faa} \mrev{with analytically (manually) computed derivatives of the flux function.}
Later, the \textit{approximate LW procedure} was introduced in~\cite{Zorio2017,Burger2017} using FD approximations.
The appropriate FD formula has to be used depending on the desired order \mrev{of} accuracy.
For fourth-order accuracy,~\cite{Zorio2017,Burger2017} used the FD approximations
\begin{equation} \label{eq:alw}
\begin{split}
\vf^{(1)} &= \frac{1}{12} [ - \pf(\vu + 2 \vu^{(1)}) + 8 \pf(\vu + \vu^{(1)}) - 8 \pf(\vu - \vu^{(1)}) + \pf(\vu - 2 \vu^{(1)}) ], \\
\vf^{(2)} &= \pf(\vu + \vu^{(1)} + \frac{1}{2} \vu^{(2)}) - 2 \pf(\vu) + \pf(\vu - \vu^{(1)} + \frac{1}{2} \vu^{(2)}), \\
\vf^{(3)} &= \frac{1}{2} \left[\pf \left(\sum_{m \subindex 0}^3 \frac{2^m}{m!}\vu^{(m)}\right)
-2 \pf \left( \sum_{m \subindex 0}^3 \frac{1}{m!}\vu^{(m)} \right)
 + 2 \pf \left( \sum_{m \subindex 0}^3 \frac{(-1)^m}{m!}\vu^{(m)}\right)
- \pf \left(\sum_{m \subindex 0}^3 \frac{(-2)^m}{m!}\vu^{(m)} \right) \right].
\end{split}
\end{equation}
As noted in~\cite{Zorio2017,Burger2017},~\eqref{eq:faa} requires the explicit computation of analytical Jacobians of the flux function which are different for different fluxes.
This is not the case with the approximate LW procedure as it only requires evaluations of the flux functions (e.g., as in~\eqref{eq:alw}).
Consequently, it was also shown in~\cite{Burger2017} that the approximate LW procedure is faster than using~\eqref{eq:faa}.
In this work, we propose to use
automatic differentiation (AD) to compute the temporal derivatives of the flux function.
\mrev{A basic introduction to AD is given in Appendix~\ref{sec:ad.appendix}.
The reader is referred to~\cite{andreas2008} for a comprehensive overview.}
\mrev{Following the discussion in Appendix~\ref{sec:ad.appendix}, for} this work, we understand AD to be a
technique to compute directional derivatives without explicitly constructing the Jacobian tensor \mrev{during the process}.
By being Jacobian free,
the AD approach \mrev{for computing the temporal derivatives of the flux}
does not suffer from the performance issues of using~\eqref{eq:faa} described in~\cite{Zorio2017,Burger2017}.
Since the first-order flux derivative can be computed as $\vf^{(1)} = \pf'(\uu) \vu^{(1)}$, it can be computed by AD as it is the derivative of $\pf$ in the direction of $\vu^{(1)}$.
For higher-order derivatives, \mrev{we use a recursive approach, that is inspired by the nested approach in~\cite{Siskind2008} and related to the recursive approach in~\cite{Szirmay-Kalos_2021}.
}
To compute the flux derivative of order \mrev{$m$}, we define \textit{derivative bundles} \mrev{which are collections of the solution and its temporal derivatives up to order $m$}:
\[
\cu^{(m)} := (\vu, \vu^{(1)}, \ldots, \vu^{(m)}), \qquad \tcu^{(m)} := (\vu^{(1)}, \ldots, \vu^{(m)}).
\]
\mrev{It can be seen from~\eqref{eq:faa} that the exact $m^\text{th}$ flux derivative is a function of the bundle $\cu^{(m)}$ allowing us to write $\vf^{(m)} = \vf^{(m)}(\cu^{(m)})$.
Since automatic differentiation can compute directional derivatives (Appendix~\ref{sec:ad.appendix}), we use the chain rule to observe that for a general function $\bb = \bb(\ba(t))$, we have $\ud_t \bb = \partial_{\ba} \bb \ud_t \ba$, showing that the temporal derivative of $\bb$ can be computed as its directional derivative with respect to $\ba$ in the direction $\ud_t \ba$.
By using this application of the chain rule to compute $\vf^{(m)}$, we obtain
\begin{equation*}
\vf^{(m)} = \Delta t \partial_t (\vf^{(m-1)}(\cu^{(m-1)})) = \mrev{\Delta t \partial_{\cu^{(m-1)}} {\vf^{(m-1)}}\partial_t \cu^{(m-1)}}.
\end{equation*}
Further, by using $\Delta t \partial_t \cu^{(m-1)} = \tcu^{(m)}$, the $m^\text{th}$ derivative of the flux can be computed as
\begin{equation} \label{eq:enzyme.ad}
\vf^{(m)} = \partial_{\cu^{(m-1)}} {\vf^{(m-1)}} \tcu^{(m)}.
\end{equation}
By~\eqref{eq:enzyme.ad}, the function $\vf^{(m)}$ is defined in terms of a directional derivative of $\vf^{(m-1)}$ with respect to the vector $\cu^{(m-1)}$ in the direction $\tcu^{(m)}$, and can thus be computed by a recursive application of AD.}
Since~\eqref{eq:enzyme.ad} gives $\vf^{(m)}$ exactly, it is equivalent to Faà di Bruno's formula~\eqref{eq:faa}.

\mrev{
\subsection{Point-wise AD}
Algorithm~\ref{alg:ad} describes the \mrev{Cauchy-Kovalevskaya procedure to compute the} \mrev{local} approximation of the time average flux~\eqref{eq:tavgflux}, where the temporal derivatives $\pf_t, \pf_{tt}, \cdots$ are computed using AD.
AD is applied to the flux function and the derivatives are thus computed point-wise.
This will be referred to as the \textit{point-wise AD} implementation.
The implementation of the approximate LW procedure is similar to Algorithm~\ref{alg:ad}, with the only difference being that FD approximations are used instead of AD.
Thus, if the \texttt{derivative\_bundle} functions in Algorithm~\ref{alg:ad} are replaced by FD approximations, the approximate Lax-Wendroff procedure will be obtained.

\mrev{
\subsubsection{Choice of library for point-wise AD}
The basic details about how AD libraries work are given in Appendix~\ref{sec:ad.appendix}.
All AD libraries can be used to compute first order directional derivatives.
Thus, by using~\eqref{eq:enzyme.ad}, any AD library can be used to compute higher-order derivatives of the flux function.
Since our code \texttt{Tenkai.jl} is a Julia package, we use AD libraries available in Julia.
The usage of an AD library like \texttt{ForwardDiff.jl}~\cite{revels2016}, which performs AD by dual arithmetic~\cite{andreas2008}, is known to lead to an exponential increase in the computational cost as the order of derivative increases~\cite{tan2022thesis}.
A more efficient alternative is Taylor arithmetic (Chapter 13 of~\cite{andreas2008}, \mrev{also discussed briefly in Appendix~\ref{sec:ad.appendix}}).
Taylor arithmetic has been implemented in \texttt{TaylorDiff.jl}~\cite{tan2022taylordiff,tan2022thesis}, in which
\mrev{the rules for high order derivatives are directly implemented (Appendix~\ref{sec:ad.appendix}).}
In this work, we found \texttt{Enzyme.jl}~\cite{moses2020instead} to be the fastest Julia library to compute the point-wise high-order derivatives of the flux function while making use of~\eqref{eq:enzyme.ad}.
\texttt{Enzyme.jl} can compute the high-order derivatives without an exponential increase in the computational cost.
This is because it computes the derivative expressions at the time of compilation, and is thus efficient without requiring high order derivative rules to be specified.
It is illustrated in Appendix~\ref{sec:ad.appendix} with an example.
As a part of this work, the function \texttt{derivative\_bundle} was developed using \texttt{Enzyme.jl} in the Julia \cite{bezanson2017julia} package \texttt{TowerOfEnzyme.jl}~\cite{towerofenzyme}, implementing~\eqref{eq:enzyme.ad} with the derivative bundle $\cu^{(m)}$.
}
}
\begin{algorithm}[!ht]
\begin{minted}[fontsize=\small,,highlightlines={10,18},highlightcolor=yellow]{julia}
for i in 1:N+1 # loop over solution points to compute ut
  f_node = flux(u[i]) # compute f_node = flux(u[i])
  F[i] = f_node
  for ii in 1:N+1 # compute ut = - divergence(f)
    ut[ii] += -D[ii, i] * f_node
  end
end
for i in 1:N+1 # loop over solution points to compute utt
  # Apply AD on flux to compute ft_node
  ft_node = derivative_bundle(flux, (u[i], ut[i]))
  F[i] += 0.5 * ft_node
  for ii in 1:N+1 # compute utt = - divergence(ft)
    utt[ii] += -D[ii, i] * ft_node
  end
end
for i in 1:N+1 # loop over solution points to compute uttt
  # Apply AD on flux to compute ftt_node
  ftt_node = derivative_bundle(flux, (u[i], ut[i], utt[i]))
  F[i] += 1/6 * ftt_node
  for ii in 1:N+1 # compute uttt = - divergence(ftt)
    uttt[ii] += -D[ii, i] * ftt_node
  end
end
...
res += D * F
\end{minted}
\caption{Local time average flux~\eqref{eq:tavgflux} computation using point-wise AD.}\label{alg:ad}
\end{algorithm}
\mrev{
\subsection{Element-wise AD}
An alternative implementation is possible by differentiating the function that computes the divergence of the flux at all solution points in an element.
This divergence function is given in Algorithm~\ref{alg:compute.flux.div}, and it consists of a loop over all solution points to compute the divergence of the flux using a derivative matrix.
This function is differentiated by using the \texttt{derivative\_bundle!} function which is developed for differentiating element-wise functions\footnote{The general difference between the \texttt{derivative\_bundle!} and the \texttt{derivative\_bundle} function is that the former allows us to differentiate functions that store their output into another data structure. This is applicable for the \texttt{compute\_fluxes!} function which stores the flux computation in another array \texttt{f}.} using~\cite{tan2022thesis,tan2022taylordiff} and has been contributed to the Julia package \texttt{Tenkai.jl}~\cite{tenkai}.
We used \texttt{TaylorDiff.jl} in this case as we found it be to the fastest library for this implementation.
However, the benefit of using Element-wise AD was found to be hardware dependent as discussed in Section~\ref{sec:2deuler}.
}

\begin{algorithm}[!ht]
\begin{minted}[fontsize=\small]{julia}
function compute_fluxes!(f, u)
  for i in 1:N+1 # loop over solution points to compute the flux
    f[i] = flux(u[i])
  end
end

function compute_div!(ut, f)
  for i in 1:N+1 # loop over solution points to compute ut
    for ii in 1:N+1 # compute ut = - divergence(f)
      ut[ii] += -D[ii, i] * f[i]
    end
  end
end
\end{minted}
\caption{\mrev{Functions to compute fluxes and their divergence at all solution points in an element. The \texttt{compute\_fluxes!} function is differentiated by AD to perform the Cauchy-Kovalevskaya procedure.}} \label{alg:compute.flux.div}
\end{algorithm}

\begin{algorithm}[!ht]
\begin{minted}[fontsize=\small]{julia}
# compute flux and use it to compute ut = -divergence(f) at all solution points
compute_fluxes!(f, u) # Stores the flux in f
compute_div!(ut, f)
# apply AD on compute_fluxes! to obtain ft at all solution points
ft = derivative_bundle!(compute_fluxes!, (u, ut), cache)
compute_div!(utt, ft) # utt = -divergence(ft)
# apply AD on compute_fluxes! to obtain ftt at all solution points
ftt = derivative_bundle!(compute_fluxes!, (u, ut, utt))
compute_div!(uttt, ftt) # uttt = -divergence(ftt)
...
F = f + 0.5 * ft + 1/6 * ftt + ...
res += D * F
\end{minted}
\caption{\mrev{Local time average flux computation using AD, element-wise version.}}\label{alg:ad.inplace}
\end{algorithm}
\mrev{\subsection{Benefits of the AD approach over the approximate LW procedure}
Using the admissibility enforcing procedure of~\cite{babbar2024admissibility,babbar2025generalized} which is also described in Section~\ref{sec:flux.correction}, we can ensure that the final evolved solution $\uu^{n+1}$ obtained by the LW scheme is admissible.
The approximate Lax-Wendroff procedure is performed by evaluating the flux function at states that are perturbed using the temporal derivatives of the solution in order to perform FD approximations~\eqref{eq:alw}.
Thus, for conservation laws whose fluxes can only be computed on admissible states, the approximate LW procedure requires admissibility of these perturbed states which is not ensured by the procedure described in Section~\ref{sec:flux.correction}.
This is also the case for several other admissibility preserving techniques for Lax-Wendroff methods that include~\cite{moe2017,Chen2016} where an admissibility preserving scheme is obtained by modifying the time averaged numerical flux.
This issue was already noted in~\cite{basak2025} where the LW schemes with the approximate Lax-Wendroff procedure were studied for relativistic hydrodynamics (RHD) equations whose fluxes can only be computed on physically admissible solutions whose velocity has a magnitude of less than 1.
A potential solution would also be to use a MOOD scheme like in~\cite{Dumbser2014,Dumbser2019} where the high order solution is discarded in \textit{troubled regions} and replaced by a more robust low order solution on subcells.
Thus, a MOOD scheme can be developed which will not perform the approximate Lax-Wendroff procedure in regions where the perturbed states are inadmissible, and it will be an interesting direction for future work.
However, the admissibility enforcing procedures of~\cite{babbar2024admissibility,babbar2025generalized,moe2017,Chen2016} are still of relevance as these schemes do not discard the high order solution completely in troubled regions, but only modify the inter-element fluxes.
The work of~\cite{basak2025} remedied this issue by performing further limiting of the perturbed states in~\eqref{eq:alw} to ensure their admissibility.
The AD approach avoids this issue entirely as the flux derivatives are only applied to the solution $\uu$ whose admissibility can be ensured by the technique of Section~\ref{sec:flux.correction} or the procedures proposed in the other references~\cite{babbar2024admissibility,babbar2025generalized,moe2017,Chen2016}.

The simplicity of the AD approach can be seen in Algorithm~\ref{alg:ad} by the fact that the same \texttt{derivative\_bundle} function is used to compute the derivatives of any order, while the approximate LW procedure requires different FD approximations for different orders~\eqref{eq:alw}.

A performance improvement of up to 61\% using AD has been observed over the approximate LW procedure in our numerical experiments (Tables~\ref{tab:wct.intel},~\ref{tab:wct.mac}).
}

\section{Numerical results} \label{sec:results}

Next, we demonstrate admissibility preservation of the AD approach and compare its performance with the approximate LW procedure.
The LW scheme with AD is implemented in \texttt{Tenkai.jl}~\cite{tenkai}.
All code required to reproduce our numerical experiments, along with the resulting data, is available in the repository~\cite{babbar2025automaticRepro}.

\subsection{Admissibility preservation} \label{sec:admissibility}
\begin{figure}
\begin{center}
\begin{tabular}{cc}
\includegraphics[width=0.45\textwidth]{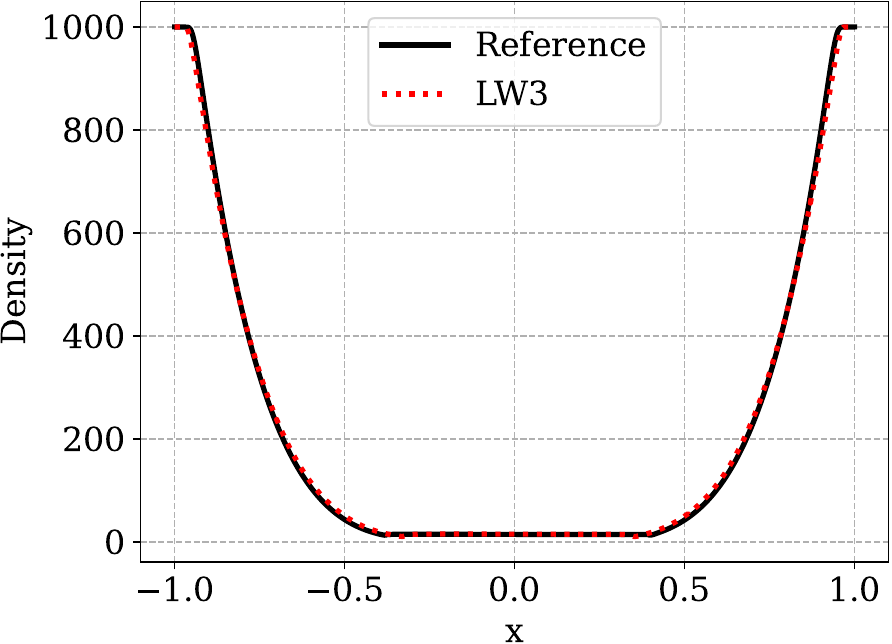} &
\includegraphics[width=0.45\textwidth]{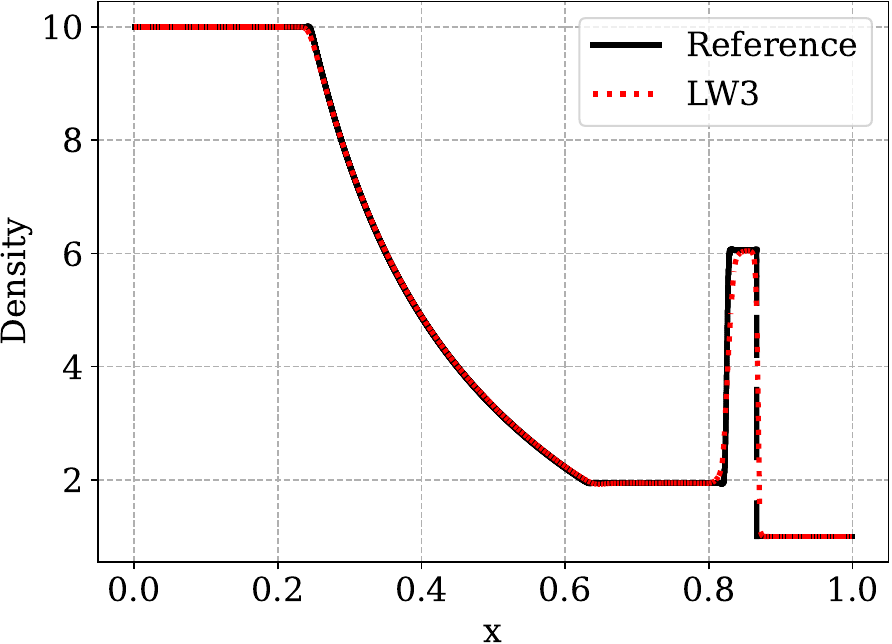} \\
(a) & (b)
\end{tabular}
\end{center}
\caption{Challenging Riemann problems for conservation laws whose fluxes can only be computed on admissible states solved using the LW-AD scheme with polynomial degree $N=3$ (LW3). (a) Double rarefaction Riemann problem for isentropic Euler equations at $t=0.2$, (b) Riemann problem for relativistic hydrodynamics from~\cite{ryu2006equation}, numbered first in~\cite{basak2025bound}.}
\label{fig:admissibility}
\end{figure}

The isentropic Euler equations are given by~\cite{Toro2009}
\mrev{\begin{equation}
\begin{pmatrix} \rho \\ \rho v \end{pmatrix}_t +\begin{pmatrix} \rho v \\ \rho v^2 + p(\rho) \end{pmatrix}_x = \bzero, \label{eq:isentropic}
\end{equation}
where the conservative variables $\rho$, $\rho v$ are the density and momentum, respectively, $v$ is the velocity, and $p(\rho) = \rho^\gamma$ is the pressure where $\gamma$ is the gas constant taken to be 1.4.
Note that the flux function requires the computation of pressure which can only be done when the density is positive, that is, the solution is physically admissible~\eqref{eq:uad.form}.}
This is different from the compressible Euler equations~\cite{Toro2009}, where the flux function can be computed even when the physical solution is inadmissible.
We solve a Riemann problem whose solution consists of two rarefaction waves. The density is initially set to $1000$ everywhere, and the velocity to $-3.9$ in the left state and $+3.9$ in the right state.
\mrev{The subcell based blending limiter~\cite{babbar2024admissibility} is used to suppress spurious oscillations near the discontinuities along with the admissibility enforcing procedure of~\cite{babbar2024admissibility}.}
The approximate LW procedure crashes as the density of one of the states in~\eqref{eq:alw} turns negative despite the admissibility correction from~\cite{babbar2024admissibility}.
However, the AD approach remains stable as it only requires the admissibility of the solution $\uu$ which is ensured through the procedure of\mrev{~\cite{babbar2024admissibility,babbar2025generalized}, also described in Section~\ref{sec:flux.correction}}. The second example is a relativistic hydrodynamics (RHD)~\cite{sokolov2001simple} problem studied for LWFR schemes initially in~\cite{basak2025bound} and later for general equations of state (EOS) in~\cite{basak2025}.
\mrev{The RHD equations are given in $d$-dimensions as~\cite{basak2025bound},
\begin{align}\label{RHD_equation}
    \frac{\partial \boldsymbol{u}}{\partial t}+\sum_{i=1}^d\frac{\partial \boldsymbol{f}_i(\boldsymbol{u})}{\partial x_i}=0,
\end{align}
where $\boldsymbol{u}$ is the vector of conserved variables and $\boldsymbol{f}_i$ is the flux vector in $x_i$ direction. Specifically, these vectors are defined by,
\begin{align}
\begin{split}
    \boldsymbol{u}&=(D, m_1,m_2,\dots,m_d, E)^\top,\\
    \boldsymbol{f}_i(\boldsymbol{u})&=(Dv_i, m_1v_i+p\delta_{1,i}, m_2v_i+p\delta_{2,i},\dots, m_dv_i+p\delta_{d,i}, m_i)^\top,\ \ \ i=1,\dots, d,
\end{split}
\end{align}
where $D=\rho\Gamma$ is the relativistic density, $\boldsymbol{m}=(m_1,m_2,\dots,m_d)^\top$ is the momentum density, $E={\rho h\Gamma^2}-p$ is the energy density where $\rho$, $p$ and $\boldsymbol{v}=(v_1, \dots, v_d)^\top$ are the rest-mass density, kinetic pressure and fluid velocity, respectively. The conservative variable $\boldsymbol{m}$ is expressed in terms of the primitive variables as $\boldsymbol{m}=\rho h \boldsymbol{v}\Gamma^2$, where $\Gamma=\frac{1}{\sqrt{1-|\boldsymbol{v}|^2}}$ is the Lorentz factor and $h$ is the specific enthalpy. The above equations need to be closed using an equation of state. Several equations of state exist for RHD equations, and for our test case, we use the one from~\cite{ryu2006equation} given by $h = \frac{2(6p^2 + 4p\rho +\rho^2)}{\rho (3p+ 2\rho)}$.
For the solution of RHD equations to be admissible, the density and pressure need to be positive and the magnitude of fluid velocity to be less than one, which is the scaled speed of light.
Since the flux function requires the computation of the Lorentz factor, similar to the isentropic Euler equations~\eqref{eq:isentropic}, the flux can only be computed on physically admissible states.}
We consider the 1-D Riemann problem from~\cite{ryu2006equation} \mrev{given by
\[
    (\rho, v_1, p) = \begin{cases}
        (10, 0, 13.3), \quad &\text{if}\ x<0.5,\\
        (1, 0, 10^{-6}), \quad &\text{if}\ x>0.5.
    \end{cases}
\]
}
In this test, the approximate LW procedure crashes as one of the states where the flux evaluation is needed in~\eqref{eq:alw} is not admissible.
In contrast, the AD approach works fine and the solution with the LW-AD scheme is shown in Figure~\ref{fig:admissibility}b.
\mrev{This test was used in~\cite{basak2025} where additional scaling was performed to keep the perturbed states in the approximate LW procedure~\eqref{eq:alw} admissible.}
The LW-AD results are compared to a first-order finite volume solution obtained on a very fine grid and show reasonable agreement.

\subsection{2-D compressible Euler equations} \label{sec:2deuler}

\begin{figure}
\begin{center}
\begin{tabular}{cc}
\includegraphics[width=0.45\textwidth]{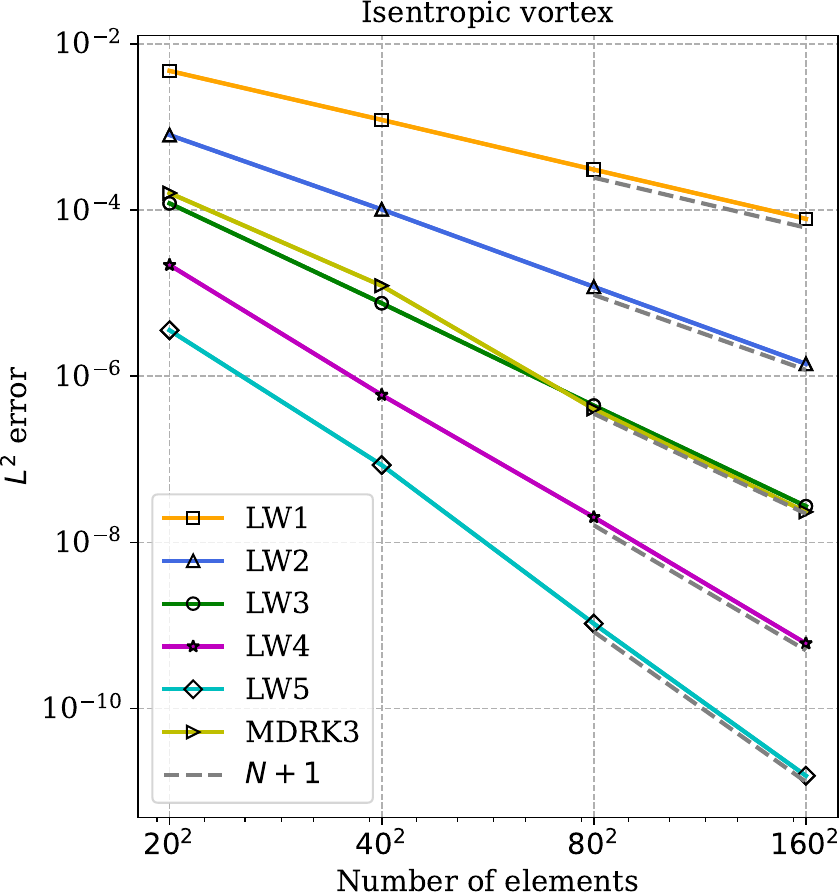} &
\includegraphics[width=0.45\textwidth]{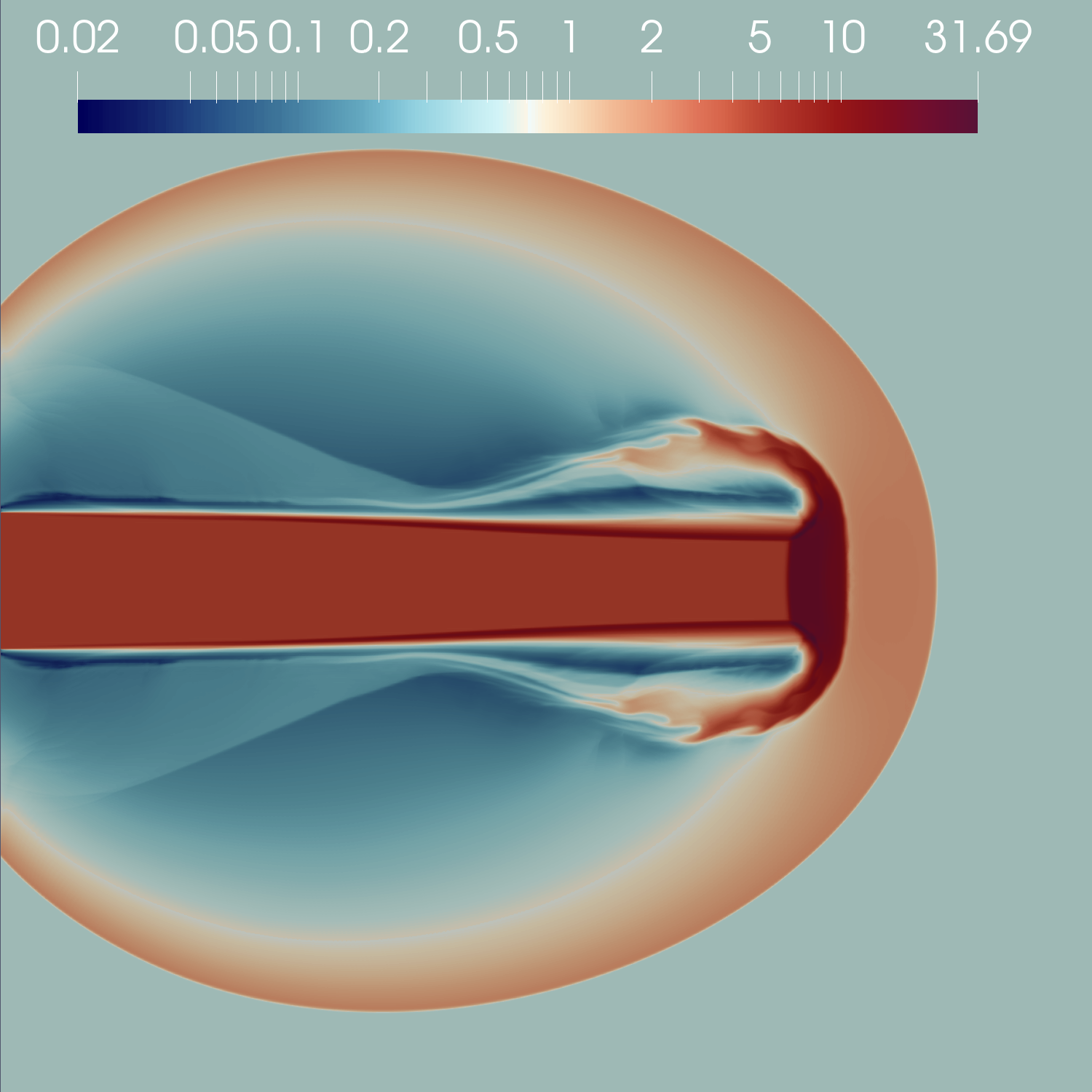} \\
(a) & (b)
\end{tabular}
\end{center}
\caption{Results with compressible Euler's equations obtained using the LW-AD scheme with different polynomial degrees $N$ (LW-$N$) and the multiderivative Runge-Kutta flux reconstruction scheme using AD using polynomial degree $N=3$ (MDRK3). (a) Convergence plot of isentropic vortex, (b) Density profile of Mach 2000 astrophysical jet. \label{fig:2deuler} }
\end{figure}

We now test the LW-AD scheme for the 2-D compressible Euler equations
\mrev{\begin{equation}\label{eq:2deuler}
\pd{}{t} \begin{pmatrix}
\rho \\
\rho v_1 \\
\rho v_2 \\
E
\end{pmatrix} +
\pd{}{x} \begin{pmatrix}
\rho v_1 \\
p + \rho v_1^2 \\
\rho v_1 v_2 \\
(E+p)v_1
   \end{pmatrix} +
\pd{}{y} \begin{pmatrix}
\rho v_2 \\
\rho v_1 v_2 \\
p + \rho v_2^2 \\
(E+p)v_2
\end{pmatrix}= 0
\end{equation}
where $\rho, p$ and $E$ denote the  density, pressure and total energy of the gas, respectively and $(v_1, v_2)$ are Cartesian components of the fluid velocity.
For a polytropic gas, an equation of state which leads to a closed system is given by $E = E(\rho, v_1, v_2, p) = \frac{p}{\gamma -1}+\frac{1}{2} \rho (v_1^2 + v_2^2)$, where $\gamma > 1$ is the adiabatic constant, that will be taken as $1.4$ in the numerical tests, which is the typical value for air.
We will consider a test case consisting of} a smooth analytical solution and another with a very strong bow shock.
The smooth test is the isentropic vortex test from~\cite{Yee1999,Spiegel2016} which is also described in~\cite{babbar2022lax,babbar2024admissibility}.
Figure~\ref{fig:2deuler} shows that the AD scheme gives optimal convergence rates for LW schemes with polynomial degrees $N=1,2,3,4,5$ and the fourth-order multiderivative Runge-Kutta FR scheme of~\cite{babbar2024multi} with polynomial degree $N=3$.
The second test is a Mach 2000 astrophysical jet flow from~\cite{ha2005, zhang2010c}, also described in~\cite{babbar2024admissibility}.
The density profile of the numerical solution is shown in Figure~\ref{fig:2deuler}b.

The AD procedure also directly applies to problems involving curvilinear grids and has been implemented in the package \texttt{TrixiLW.jl}~\cite{babbar2024trixilw} using \texttt{Trixi.jl}~\cite{schlottkelakemper2021purely,schlottkelakemper2020trixi,ranocha2022adaptive} as a library.
We were these able to obtain similar results as in~\cite{babbar2025} using AD.

We use the isentropic vortex setup for a wall-clock time (WCT) comparison between the approximate LW procedure (\textbf{FD},~\eqref{eq:alw}) and the AD approach (\textbf{AD}, Algorithm~\ref{alg:ad}) on an Intel Xeon Gold 5320 CPU (Table~\ref{tab:wct.intel}) \mrev{and an Apple M3 CPU (Table~\ref{tab:wct.mac}).}
The WCT samples are taken by running the simulation thrice and taking the minimum to reduce noise.
\mrev{
With both the CPUs, the results show that the AD approach is faster than the approximate LW procedure for all polynomial degrees and mesh sizes tested.
For the intel CPU (Table~\ref{tab:wct.intel}), on the finest mesh, a speedup of $1.28, 1.15, 1.18, 1.13$ is observed for polynomial degrees $N=1,2,3,4$, respectively, using point-wise AD and a speedup of $1.32, 1.28, 1.17, 1.04$ is observed using element-wise AD.
For the Apple M3 CPU (Table~\ref{tab:wct.mac}), on the finest mesh, a speedup of $1.33, 1.07, 1.17, 1.11$ is observed for polynomial degrees $N=1,2,3,4$, respectively, using point-wise AD and a speedup of $1.61, 1.51, 1.41, 1.15$ is observed using element-wise AD.
Thus, the element-wise AD implementation is clearly better on the Apple M3 CPU, while it is comparable or slightly worse than the point-wise AD implementation on the Intel CPU.
Since both of the architectures are very different, it is beyond the scope of this work to explain this behavior.
}

\section{Conclusions} \label{sec:conclusions}

\mrev{Automatic differentiation (AD) to perform the Cauchy-Kovalevskaya step in high order Lax-Wendroff schemes has been proposed.
AD is performed using point-wise and element-wise AD.
The former is a scheme obtained simply by replacing the finite difference approximations in the approximate Lax-Wendroff procedure with AD.
The latter was obtained by differentiating the element wise flux computation function, and was seen to be faster than the point-wise version on some hardware architectures.
}
The numerical results show that the LW scheme with AD is more robust than the approximate LW procedure, without losing out on performance.
The AD approach is also simpler to implement because it does not require different finite difference approximations for different orders of accuracy, as we demonstrated by extending the LW method to the sixth order of accuracy.
The following questions are open for future research.
\mrev{Since the Lax-Wendroff method with Cauchy-Kovalevskaya procedure has also been used with other spatial discretizations like WENO~\cite{Qiu2003}, it will be interesting to see how AD performs when combined with other spatial discretizations.}
LW with AD needs to be tested for other equations including different conservation laws, but also equations with source terms~\cite{babbar2025generalized} and parabolic terms~\cite{babbar2025parabolic}.
When computing higher-order derivatives, some computations from lower-order derivatives can be reused, although this is not done in the current implementation (Algorithm~\ref{alg:ad}).
A hard-coded implementation of the LW scheme using the chain rule could reuse these intermediate results to avoid redundant operations.
However, it is not clear whether this potential optimization is possible for the proposed LW-AD scheme\mrev{s}.

\begin{table}
\centering
\begin{tabular}{clllrrrrrrrrrr}
\hline
\multicolumn{1}{|c|}{\nel} & \multicolumn{3}{|c|}{$N=1$} & \multicolumn{3}{|c|}{$N=2$} & \multicolumn{3}{|c|}{$N=3$} & \multicolumn{3}{|c|}{$N=4$} \\
\hline
\textbf{ } & \textbf{FD} & \textbf{ADP}  & \textbf{ADE} & \textbf{FD} & \textbf{ADP} & \textbf{ADE} & \textbf{FD} & \textbf{ADP} & \textbf{ADE} & \textbf{FD} & \textbf{ADP} & \textbf{ADE} \\\hline
$20^2$ & 0.009 & 0.007 & 0.007 & 0.03 & 0.02 & 0.02 & 0.1 & 0.1 & 0.1 & 0.3 & 0.2 & 0.3 \\
$40^2$ & 0.07 & 0.05 & 0.06 & 0.26 & 0.22 & 0.20 & 1.0 & 0.8 & 0.8 & 2.4 & 2.1 & 2.3 \\
$80^2$ & 0.61 & 0.48 & 0.46 & 2.11 & 1.82 & 1.63 & 7.5 & 6.4 & 6.5 & 19.3 & 16.8 & 18.9 \\
$160^2$ & 4.93 & 3.83 & 3.73 & 16.80 & 14.50 & 13.10 & 61.2 & 51.8 & 52.2 & 158.0 & 139.0 & 151.0 \\\hline
$\frac{\textbf{FD}_{160^2}}{\textbf{x}_{160^2}}$ & 1.0 & 1.28 & 1.32 & 1.0 & 1.15 & 1.28 & 1.0 & 1.18 & 1.17 & 1.0 & 1.13 & 1.04
\end{tabular}
\caption{Wall clock time performance on an Intel CPU (in seconds). $\nel$: Number of elements, $N$: Polynomial degree, FD: approximate Lax-Wendroff procedure~\eqref{eq:alw}, ADP: point-wise AD (Algorithm~\ref{alg:ad}), ADE: element-wise AD (Algorithm~\ref{alg:ad.inplace}).
The last row shows the ratio of the WCT of the approximate LW procedure to the AD approach on the finest mesh, measuring the speedup of the AD approach.}
\label{tab:wct.intel}
\end{table}

\begin{table}
\centering
\begin{tabular}{clllrrrrrrrrrr}
\hline
\multicolumn{1}{|c|}{\nel} & \multicolumn{3}{|c|}{$N=1$} & \multicolumn{3}{|c|}{$N=2$} & \multicolumn{3}{|c|}{$N=3$} & \multicolumn{3}{|c|}{$N=4$} \\
\hline
\textbf{ } & \textbf{FD} & \textbf{ADP} & \textbf{ADE} & \textbf{FD} & \textbf{ADP} & \textbf{ADE} & \textbf{FD} & \textbf{ADP} & \textbf{ADE} & \textbf{FD} & \textbf{ADP} & \textbf{ADE} \\\hline
$20^2$ & 0.005 & 0.004 & 0.004 & 0.01 & 0.01 & 0.01 & 0.06 & 0.05 & 0.04 & 0.14 & 0.13 & 0.12 \\
$40^2$ & 0.05 & 0.03 & 0.02 & 0.14 & 0.13 & 0.09 & 0.49 & 0.41 & 0.34 & 1.22 & 1.10 & 1.05 \\
$80^2$ & 0.35 & 0.27 & 0.22 & 1.17 & 1.09 & 0.76 & 3.87 & 3.26 & 2.71 & 9.62 & 8.70 & 8.32 \\
$160^2$ & 2.81 & 2.10 & 1.74 & 9.14 & 8.47 & 6.04 & 30.10 & 25.60 & 21.20 & 76.60 & 69.00 & 66.30 \\\hline
$\frac{\textbf{FD}_{160^2}}{\textbf{x}_{160^2}}$  & 1.0 & 1.33 & 1.61 & 1.0 & 1.07 & 1.51 & 1.0 & 1.17 & 1.41 & 1.0 & 1.11 & 1.15
\end{tabular}
\caption{Wall clock time performance on an Apple M3 CPU (in seconds). $\nel$: Number of elements, $N$: Polynomial degree, FD: approximate Lax-Wendroff procedure~\eqref{eq:alw}, ADP: point-wise AD (Algorithm~\ref{alg:ad}), ADE: element-wise AD (Algorithm~\ref{alg:ad.inplace}).
The last row shows the ratio of the WCT of the approximate LW procedure to the AD approach on the finest mesh, measuring the speedup of the AD approach.}
\label{tab:wct.mac}
\end{table}


\section*{Acknowledgments}

\input{funding}

\appendix

\section{\mrev{Basic implementation of automatic differentiation}} \label{sec:ad.appendix}

The idea of automatic differentiation is to teach the differentiation of some fundamental functions like the trigonometric functions, exponential and logarithmic functions, etc. along with basic rules of derivatives of functions given by
\[
(f+g)' = f' + g', \quad (f g)' = f' g + f g', \quad \left(\frac{f}{g}\right)' = \frac{f' g - f g'}{g^2}, \quad (f \circ g)' = (f' \circ g) g',
\]
to the computer so that it can compute derivatives of complicated functions built up from these basic operations.
In this work, we use forward mode AD where the flow of derivative computations follows the flow of function evaluations.
The alternative is reverse mode AD which is more efficient for functions with a large input dimension and small output dimension.
The reader can read more about forward and reverse mode AD in Chapter 3 of~\cite{andreas2008}.
A natural way to express the forward mode AD is through the use of dual arithmetic.
A general dual number is written as $x + \epsilon y$ where $\epsilon$ is a symbol formally defined to satisfy $\epsilon^2 = 0$. This is similar to ignoring the second order terms in a Taylor expansion.
By this property of $\epsilon$, addition and multiplication of dual numbers can be defined as\footnote{Since the product of two non-zero dual numbers can be zero, the dual numbers do not form a field. However, they do form a commutative algebra of dimension two over the reals.}
\begin{equation}
(x_1 + \epsilon y_1) + (x_2 + \epsilon y_2) = (x_1 + x_2) + \epsilon (y_1 + y_2), \quad
(x_1 + \epsilon y_1)(x_2 + \epsilon y_2) = x_1 x_2 + \epsilon (x_1 y_2 + x_2 y_1).
\label{eq:dual.arithmetic}
\end{equation}
If dual numbers are of the form $f(a) + \epsilon f'(a)$ and $g(a) + \epsilon g'(a)$, then the operations in~\eqref{eq:dual.arithmetic} correspond to the addition and product rules of derivatives.
For showing how chain rule is performed using dual numbers and to give a complete description, we show how dual arithmetic can be implemented through in a Julia code\footnote{We warn the reader that extending this approach to compute high order derivatives additionally requires solving the problem of perturbation confusion~\cite{manzyuk2019}.
It is a well-studied problem in the literature and solutions have already been implemented in reputed AD libraries like \texttt{ForwardDiff.jl}, allowing the usage of~\eqref{eq:enzyme.ad} with \texttt{ForwardDiff.jl}.}.
The first step is to define a new data structure \texttt{Dual} that can hold the value and derivative of a function at some fixed point.
\begin{minted}[fontsize=\small]{julia}
# Define a new data structure Dual, that inherits properties from the
# abstract type Number. It has two subfields (value and derivative)
# as double-precision floating point numbers.
struct Dual <: Number
  value::Float64
  derivative::Float64
end
\end{minted}
A \texttt{Dual} number thus contains $[v_f,d_f]$ where $v_f$ holds the \texttt{value} and $d_f$ the \texttt{derivative} of the function $f$ at some fixed point.
The basic rules of arithmetic operations on \texttt{Dual} numbers are then prescribed by giving a definition of these arithmetic operations for the \texttt{Dual} numbers.
For example, to emulate the product rule for derivatives, we consider another \texttt{Dual} number $[v_g,d_g]$ corresponding to a function $g$ at the same fixed point.
We then define multiplication of the two dual numbers following~\eqref{eq:dual.arithmetic} as $[v_f, d_f] [v_g, d_g] = [v_f v_g, v_f d_g + v_g d_f]$.
Thus, the product gives another \texttt{Dual} number where the \texttt{value} field contains the values of $fg$, and the \texttt{derivative} field contains the derivative of the product of $fg$ at the same fixed point.
The other arithmetic operations can be defined similarly.
In order to give a complete demonstration, we give a working code defining these operations for \texttt{Dual} numbers using Julia's multiple dispatch.
\begin{minted}[fontsize=\small]{julia}
import Base: +, -, *, / # Extending operations requires importing them from Base
function +(x::Dual, y::Dual) # Define addition of Dual numbers
  return Dual(x.value + y.value, x.derivative + y.derivative)
end
function -(x::Dual, y::Dual) # Define subtraction of Dual numbers
  return Dual(x.value - y.value, x.derivative - y.derivative)
end
function *(x::Dual, y::Dual) # Define multiplication of Dual numbers
  return Dual(x.value * y.value, x.value * y.derivative + x.derivative * y.value)
end
function /(x::Dual, y::Dual) # Define division of Dual numbers
  der = (x.derivative * y.value - x.value * y.derivative) / y.value^2
  return Dual(x.value / y.value, der)
end
\end{minted}
The chain rule operation is emulated using \texttt{Dual} numbers by defining applications of certain functions on the \texttt{Dual} numbers.
In general, we define the application of a function $g$ on the \texttt{Dual} number $[v_f, d_f]$ as $g([v_f, d_f]) = [g(v_f), g'(v_f) d_f]$ so that its \texttt{value} contains the value of $g \circ f$ and its derivative contains the derivative of $g \circ f$ at some fixed value.
The chain rule is thus implemented by specifying the derivatives of some fundamental functions.
\begin{minted}[fontsize=\small]{julia}
Base.convert(::Type{Dual}, x::Real) = Dual(x, 0.0) # Conversion from Real to Dual
Base.promote_rule(::Type{Dual}, ::Type{<:Real}) = Dual # Promotion for mixed operations
function Base.sin(x::Dual) # Define sine function for Dual numbers
  return Dual(sin(x.value), cos(x.value) * x.derivative) # sin(f(x)), cos(f(x)) df/dx
end
function Base.cos(x::Dual) # Define cosine function for Dual numbers
  return Dual(cos(x.value), -sin(x.value) * x.derivative) # cos(f(x)), -sin(f(x)) df/dx
end
function Base.log(x::Dual) # Define logarithm function for Dual numbers
  return Dual(log(x.value), x.derivative / x.value) # log(f(x)), df/dx * 1/f(x)
end
function Base.exp(x::Dual) # Define exponential function for Dual numbers
  e = exp(x.value)
  return Dual(e, e * x.derivative) # exp(f(x)), exp(f(x)) df/dx
end
\end{minted}
With these definitions, we compute the derivative of a function $f$ at a point $x$ in the direction $v$ as follows:
\begin{minted}[fontsize=\small]{julia}
differentiate(f, x, v) = f(Dual(x, v)).derivative # Define the differentiation function
f(x) = sin(exp(x^2)) + log(x^4 + cos(x^3)^2) # Example function
differentiate(f, 1, 2) # Differentiate at x = 1 in the direction v = 2 giving df/dx(1) * 2
\end{minted}
The above description applies to derivatives of single variable functions, but the same idea can be extended to multivariate functions by defining a dual number type that can hold vectors for values and derivatives.
In addition, high order derivatives are computed by defining a Taylor bundle type that can hold values of derivatives up to a certain order and defining rules for arithmetic operations on these Taylor bundles.
Here we simply mention the arithmetic rules that need to be defined for high order derivatives.
\[
(f+g)^{(m)} = f^{(m)} + g^{(m)}, \quad (f g)^{(m)} = \sum_{k=0}^m \binom{m}{k} f^{(k)} g^{(m-k)}, \quad \left(\frac{f}{g}\right)^{(m)} = \frac{1}{g} \left( f^{(m)} - \sum_{k=1}^{m-1} \binom{m}{k} g^{(m-k)} \left(\frac{f}{g}\right)^{(k)} \right).
\]
The exact chain rule given by Faà di Bruno's formula~\eqref{eq:faa} is not prescribed but rather the following formula is inductively taught for the fundamental functions
\[
f(g(x))^{(n)} = (f'(g(x)) g'(x))^{(n-1)} = \sum_{i=0}^{n-1} \binom{n-1}{i} (f'(g(x)))^{(i)} g^{(n-1)}(x).
\]
This idea of Taylor arithmetic is implemented in \texttt{TaylorDiff.jl}~\cite{tan2022taylordiff}, with the implementation details provided in~\cite{tan2022thesis}.
Detailed descriptions and the theory of Taylor arithmetic can be found in Chapter 13 of~\cite{andreas2008}.

\begin{figure}
\begin{center}
\begin{tabular}{cc}
\includegraphics[width=0.45\textwidth]{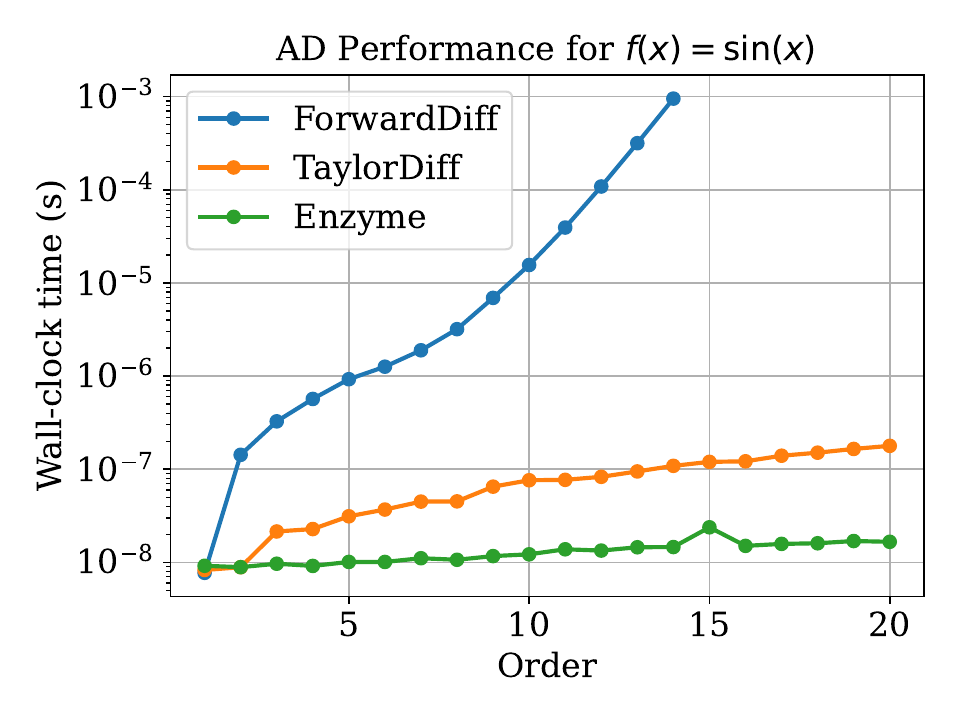} &
\includegraphics[width=0.45\textwidth]{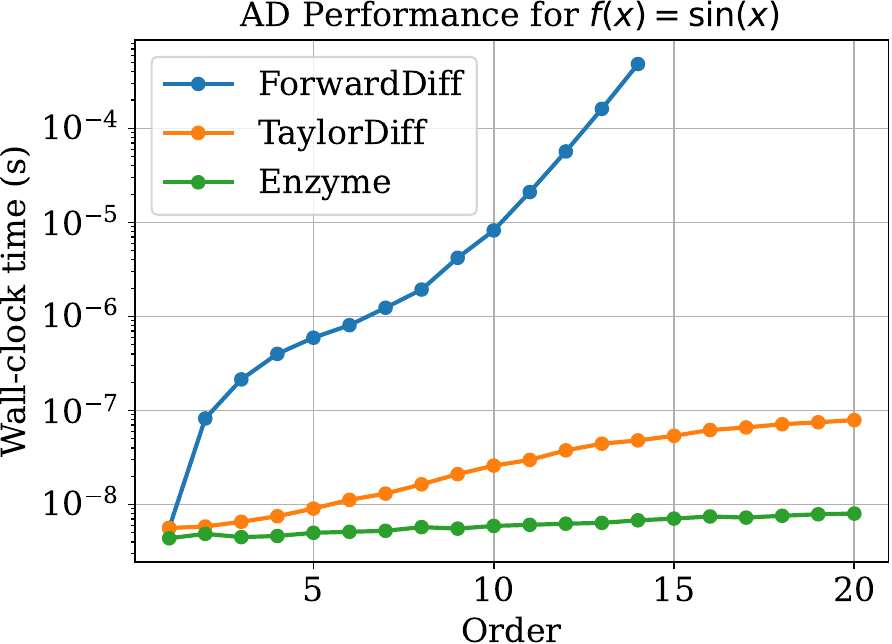} \\
(a) & (b)
\end{tabular}
\end{center}
\caption{Wall-clock time taken to compute high order derivatives of $\sin x$ at $x=1.0$ using \texttt{ForwardDiff.jl}, \texttt{TaylorDiff.jl}, and \texttt{Enzyme.jl} (a) Intel CPU, (b) Apple M3 CPU.}
\label{fig:order.vs.wct}
\end{figure}

Lastly, we discuss \texttt{Enzyme.jl} which is an AD library that works at the level of the LLVM compiler that the Julia language uses.
Enzyme is faster than the above libraries because of the following reasons.
Unlike the dual and Taylor arithmetic, Enzyme does not require any additional storage for holding the derivatives.
This avoids the potential overhead of memory allocations which can be significant when computing high order derivatives.
Enzyme interacts with LLVM code, which is the low level compiler in which Julia is written.
Following~\eqref{eq:enzyme.ad}, high order derivatives are defined recursively, and Enzyme converts each of the recursive stages to low level LLVM code allowing optimization for each derivative.
We illustrate this by showing wall clock time and some low level code obtained by high order derivative functions for the simple function $f(x) = \sin x$.
We compute high order derivatives by using differentiation functions from \texttt{ForwardDiff.jl}, \texttt{TaylorDiff.jl}, and \texttt{Enzyme.jl} and look at the generated low level LLVM code.
A perfect result in the machine code would comprise only the $\sin$ or $\cos$ function for the respective order.
We show part of the low level LLVM obtained from \texttt{Enzyme.jl} for the $20^\text{th}$ derivative of $\sin x$, hiding the boilerplate code that does not contain floating point operations.
\begin{minted}[
  fontsize=\small,
  breaklines
]{llvm}
define double @julia_nth_derivative_10736(double %0) #0 {
  ; Setting safepoints, and other semantics/boilerplate
  ; code required by Julia which does not contain any floating point operations
  %43 = call fast double @llvm.sin.f64(double %0)
  ret double %43
}
\end{minted}
We can see that the machine code contains only one \texttt{sin} call. The other libraries either fail to work at earlier orders (\texttt{ForwardDiff.jl}) or their machine code consists of more floating point operations than only the application of the sine function (\texttt{TaylorDiff.jl}).
The low level codes for the other libraries are not shown, but have been put as the reproducible Pluto.jl notebook \texttt{machine\_code.jl} in our reproducibility repository~\cite{babbar2025automaticRepro}.
In Figure~\ref{fig:order.vs.wct}, we show a log scaled plot of the wall-clock time taken to compute high order derivatives of $\sin x$ at $x=1.0$ using \texttt{ForwardDiff.jl}, \texttt{TaylorDiff.jl}, and \texttt{Enzyme.jl}.
An exponential increase in the wall clock time is seen for \texttt{ForwardDiff.jl}. Although such an increase in the computational cost is not seen for \texttt{TaylorDiff.jl}, \texttt{Enzyme.jl} is a still significantly faster.

\printbibliography

\end{document}

%% file: macros.tex
\usepackage{xparse}

\let\epsilon\varepsilon
\let\phi\varphi
\let\rho\varrho

\newcommand{\nel}{\ensuremath{N_e}}
\usepackage{xspace}

\newcommand{\emh}{e - \frac{1}{2}}
\newcommand{\eph}{e + \frac{1}{2}}
\providecommand{\subindex}{=}
\newcommand{\poly}{\mathbb{P}}
\newcommand{\uu}{\boldsymbol{u}}

\newcommand{\uep}{\uu_{e, p}}

\newcommand{\pf}{\boldsymbol{f}}
\newcommand{\F}{\boldsymbol{F}}
\newcommand{\bzero}{\boldsymbol{0}}
\newcommand{\myvector}[1]{\text{{\hspace{0.1em}}#1{\hspace{0.15em}}}}

\newcommand{\vu}{\myvector{u}}
\newcommand{\cu}{\mathcal{B}_\vu}
\newcommand{\tcu}{\tilde{\mathcal{B}}_\vu}
\newcommand{\vf}{\myvector{f}}

\newcommand{\Uad}{\mathcal{U}_{\textrm{ad}}}
\newcommand{\cdummy}{\cdot}
\newcommand{\avg}[1]{\overline{#1}}
\newcommand{\au}{\avg{\uu}}

\newcommand{\ueppone}{\uu_{e, p + 1}}

\newcommand{\uez}{\uu_{e, 0}}

\newcommand{\uepoz}{\uu_{e + 1, 0}}
\newcommand{\ueN}{\uu_{e, N}}

\newcommand{\pph}{p + \frac{1}{2}}
\newcommand{\pmh}{p - \frac{1}{2}}

\newcommand{\half}{\frac{1}{2}}

\newcommand{\Nmh}{N - \frac{1}{2}}

\usepackage[T3,T1]{fontenc}
\DeclareSymbolFont{tipa}{T3}{cmr}{m}{n}
\DeclareMathAccent{\invbreve}{\mathalpha}{tipa}{16}
\newcommand{\atu}{\invbreve{\boldsymbol{u}}}
\newcommand{\utilow}{\invbreve{\boldsymbol{u}}^{\text{low}, n + 1}}
\newcommand{\ad}{P}
\newcommand{\re}{\mathbb{R}}

\newcommand{\mrev}[1]{{\color{black}#1}}

\newcommand{\ud}{\textrm{d}}

\newcommand{\ba}{\boldsymbol{a}}
\newcommand{\bb}{\boldsymbol{b}}


%% file: abstract.tex
Lax-Wendroff methods combined with discontinuous Galerkin/flux reconstruction spatial discretization provide a high-order, single-stage, quadrature-free method for solving hyperbolic conservation laws.
In this work, we introduce automatic differentiation (AD) \mrev{for performing the Cauchy-Kowalewski procedure used} in the element-local time average flux computation step (the predictor step) of Lax-Wendroff methods.
The application of AD is similar for methods of any order and does not need positivity corrections during the predictor step.
This contrasts with the approximate Lax-Wendroff procedure, which requires different finite difference formulas for different orders of the method and positivity corrections in the predictor step for fluxes that can only be computed on admissible states.
The method is Jacobian-free and problem-independent, allowing direct application to any physical flux function.
Numerical experiments demonstrate the order and positivity preservation of the method.
Additionally, performance comparisons indicate that the wall-clock time of automatic differentiation is always on par with the approximate Lax-Wendroff method.

%% file: funding.tex
AB, VC, and HR were supported by the Deutsche Forschungsgemeinschaft
(DFG, German Research Foundation, project number 528753982
as well as within the DFG priority program SPP~2410 with project number 526031774)
and the Daimler und Benz Stiftung (Daimler and Benz foundation,
project number 32-10/22).
AB was also supported by the Alexander von Humboldt Foundation.
MSL was supported by the Deutsche Forschungsgemeinschaft
(DFG, German Research Foundation, project number 528753982
as well as within the DFG research unit FOR~5409 "SNuBIC" with project number 463312734)
and the BMBF project "ADAPTEX".
We thank Sujoy Basak for the permission to use his RHD code, and for comments on the manuscript.